\documentclass{article}

\usepackage{amsmath}
\usepackage{graphicx}
\usepackage{amssymb}
\usepackage{here}
\usepackage[sort]{natbib}
\usepackage{apalike}

\newtheorem{theo}{Theorem}

\def\qed{\hfill $\Box$}
\newcommand{\ten}[1]{#1^\top}
\newcommand{\I}[1]{\mathrm{I}_{#1}}
\newcommand{\argmin}{\mathop{\rm arg~min}\limits}
\newcommand{\inv}[1]{#1^{-1}}

\newcommand{\aic}{\mathrm{AIC}}
\newcommand{\aicc}{\mathrm{AICc}}

\newcommand{\zklic}{\mathrm{ZKLIC}}

\newcommand{\mcp}{\mathrm{MCp}}
\newcommand{\zmcp}{\mathrm{ZMCp}}

\newcommand{\tjbay}{\hat{\Theta}_{J,\mathrm{B}}}

\newcommand{\bjbay}{\hat{B}_{J,\mathrm{B}}}
\newcommand{\bjmle}{\hat{B}_J}

\newcommand{\E}[1]{\mathrm{E}\left[#1\right]}
\newcommand{\Es}[2]{\mathrm{E}_{#1}\left[#2\right]}

\newcommand{\srisk}[3]{R_{\mathrm{S}}(#1,#2,#3)}
\newcommand{\klrisk}[3]{R_{\mathrm{KL}}(#1,#2,#3)}
\newcommand{\rank}{\mathrm{rank}}

\newcommand{\invjxisx}{(\ten{x}_{J,i}\inv{S}_Jx_{J,i})^{-1}}

\newcommand{\wjicsjwjic}{\ten{w}_{J,i}\inv{S}_Jw_{J,i}}
\newcommand{\wjicsfwjic}{\ten{w}_{J,i}\inv{S}_Fw_{J,i}}
\newcommand{\invwjicsfwjic}{(\ten{w}_{J,i}\inv{S}_Fw_{J,i})^{-1}}
\newcommand{\invwjicsjwjic}{(\ten{w}_{J,i}\inv{S}_Jw_{J,i})^{-1}}

\newcommand{\kj}{k_J}
\newcommand{\kf}{k_F}
\newcommand{\kjt}{k_{J*}}
\newcommand{\Kj}{K_J}
\newcommand{\Pj}{P_J}
\newcommand{\Pjt}{P_{J_*}}
\newcommand{\cji}{c_{J,i}}

\newcommand{\Cj}{C_J}

\newcommand{\dji}{d_{J,i}}

\newcommand{\Dj}{D_J}

\newcommand{\lji}{\lambda_{J,i}}

\newcommand{\Lj}{\Lambda_{J}}
\newcommand{\eji}{\eta_{J,i}}

\newcommand{\xj}{X_J}
\newcommand{\xjic}{x_{J,i}}

\newcommand{\tj}{\Theta_J}
\newcommand{\tjt}{\Theta_{J_*}}
\newcommand{\tjic}{\Theta_{J,i\cdot}}
\newcommand{\qj}{Q_J}
\newcommand{\aj}{A_J}
\newcommand{\af}{A_F}
\newcommand{\ajt}{A_{J_*}}
\newcommand{\wj}{W_J}
\newcommand{\wji}{w_{J,i}}
\newcommand{\wjic}{w_{J,i}}
\newcommand{\bj}{B_J}

\newcommand{\bjt}{B_{J_*}}
\newcommand{\zj}{Z_J}
\newcommand{\zf}{Z_F}
\newcommand{\sj}{S_J}
\newcommand{\lj}{L_J}
\newcommand{\ljt}{L_{J_*}}
\newcommand{\mj}{M_J}
\newcommand{\mjt}{M_{J_*}}
\newcommand{\isf}{\inv{S_F}}
\newcommand{\isj}{\inv{S_J}}

\newcommand{\Vji}[2]{V_{J,i,#1}^{#2}}
\newcommand{\Rji}[1]{R_{J,i}^{(#1)}}

\newcommand{\vji}[1]{v_{J,i,#1}}

\newcommand{\vjivji}[2]{\ten{v_{J,i,#1}}v_{J,i,#2}}

\newcommand{\uji}[1]{u_{J,i,#1}}
\newcommand{\Sigmajmle}{\hat{\Sigma}_{J}}
\newcommand{\half}{\frac{1}{2}}

\newcommand{\rj}{r_J}

\def\m{\quad}
\def\n{\enskip}

\def\diag{\mathop{\rm diag}\nolimits}
\def\tr{\mathop{\rm tr}\nolimits}

\def\sqr#1#2{{\vcenter{\hrule height.#2pt
    \hbox{\vrule width.#2pt height#1pt \kern#1pt \vrule width.#2pt}
    \hrule height.#2pt}}}

\title{Generalized ridge estimator and model selection criterion 
in multivariate linear regression}
\author{Yuichi Mori\footnote{Department of Mathematical 
and Computing Sciences Graduate School of Information 
Science and Engineering Tokyo Institute of Technology}, 
Taiji Suzuki$^*
\footnote{PRESTO,JST}$
}
\date{}

\begin{document}
\bibliographystyle{apalike}

\maketitle
\renewcommand{\abstractname}{Abstract}
\begin{abstract}
We propose new model selection criteria based on generalized 
ridge estimators dominating the maximum likelihood estimator 
under the squared risk and the Kullback-Leibler risk 
in multivariate linear regression. Our model selection criteria have the 
following favorite properties: consistency, unbiasedness, uniformly minimum variance. 
Consistency is proven under an asymptotic structure 
$\frac{p}{n}\to c$ where $n$ is the sample size and $p$ is the parameter dimension 
of the response variables.
In particular, our proposed class of estimators dominates the maximum likelihood estimator 
under the squared risk even when the model does not include the true model. 
Experimental results show that the risks of  
our model selection criteria are smaller than the ones based on the maximum likelihood 
estimator and that our proposed criteria specify the true model under some conditions.
\end{abstract}

\section{Introduction}
\m\n Model selection criteria such as AIC (\cite{akaike1971information}) 
and Cp (\cite{mallows1973some}) have been used in various applications and 
their theoretical properties have been extensively studied. 
We consider a model selection problem in a multivariate linear regression 
based on a kind of generalized ridge estimators.
The multivariate linear regression considered in this paper has $p$ response variables on 
a subset of $k$ explanatory variables, and  
the response is contaminated by a multivariate normal noise. 
This model, in which the response is multivariate, is thus an extension of multiple linear regression, where the response is univariate.
Applications of multivariate linear regression include genetic data analysis, (e.g., 
\cite{gharagheizi2008qspr}) and multiple brain scans (e.g., \cite{basser1998simplified}).
Multivariate linear regression is written as
\begin{eqnarray*}
Y\sim {\cal N}_{n\times p}(AB,\Sigma\otimes\I{n}),
\end{eqnarray*}
where $Y$ is an $n\times p$ observation matrix of $p$ response variables, 
$A$ is an $n\times k$ observation matrix of $k$ explanatory variables, 
 $B$ is a $k\times p$ unknown matrix of regression coefficients, 
$\Sigma$ is a $p\times p$ unknown covariance matrix,
$k$ is a nonstochastic number, and $n$ is the sample size. We assume that, 
for all  $n\geq k$, $n-p-k-1>0$ and that $\rank(A)=k$.

The purpose of the model selection problem is to select an appropriate 
subset of regression coefficients. 
Suppose that $J$ denotes a subset of the index set of coefficients $F=\{1,...,k\}$.
${\cal J}$ denotes the power set of $F$, and $\kj$ denotes 
the number of elements that $J$ contains, that is, $\kj=|J|$. 
Then, the candidate model corresponding to 
the subset $J$ can be expressed as
\begin{eqnarray*}
Y\sim {\cal N}_{n\times p}(\aj\bj,\Sigma\otimes\I{n}),
\end{eqnarray*}
where $\aj$ is an $n\times \kj$ matrix consisting of the columns 
of $A$ indexed by the elements of $J$, and $\bj$ is a $\kj\times p$ 
unknown matrix of regression coefficients. We assume that 
the candidate model corresponding to $J_*\in{\cal J}$ is the true model.

One way to perform model selection in multivariate linear regression is to apply the well known model selection criteria such as 
AIC (\cite{akaike1971information}), AICc (\cite{bedrick1994model}), 
Cp (\cite{mallows1973some}), and MCp (\cite{fujikoshi1997modified}).
These criteria are unbiased or asymptotic unbiased estimators 
of the squared risk and the Kullback Leibler risk that are defined as follows:
\begin{eqnarray*}
\srisk{B}{\Sigma}{\Phi} &=&
\E{\tr\left(\inv{\Sigma}\ten{(B-\Phi)}\ten{A}A(B-\Phi)\right)},\\
\klrisk{B}{\Sigma}{\hat{f}} &=& \Es{\tilde{Y},Y}{
\log\left(
\frac{f(\tilde{Y}|B,\Sigma)}{\hat{f}(\tilde{Y}|Y)}
\right)},
\end{eqnarray*}
where $\Phi$ is an estimator of $B$, $f$ is the true probability density of $Y$, and 
$\hat{f}$ is a predictive density of $\tilde{Y}$ conditional to $Y$ 
where $\tilde{Y}$ is the independent copy of $Y$. 
Cp and MCp are unbiased estimators of the squared risk of the 
maximum likelihood estimator, and 
AIC and AICc are asymptotic unbiased and unbiased estimators, respectively, of 
the Kullback-Leibler risk of the maximum likelihood estimator. 
In particular, it is shown that MCp and AICc are uniformly minimum variance unbiased 
estimators of their corresponding risks (\cite{davies2006estimation}).
One important property of a model selection criterion is consistency, that is, 
it selects the true model asymptotically in probability. 
\cite{fujikoshi2014consistency} showed consistency of AIC, AICc , Cp and MCp 
under an asymptotic structure $\frac{p}{n}\rightarrow c$\ and some conditions 
in multivariate linear regression, although they are not consistent 
in usual univariate settings. 

Besides model selection criteria corresponding to the maximum likelihood estimator as 
introduced above, some authors have studied those corresponding to other 
estimators that might dominate the maximum likelihood estimator. 
\cite{yanagihara2010unbiased} investigated an unbiased 
estimator of the squared risk of the ridge estimator. 
They developed a model selection criterion to select the 
model candidate and the parameter of the ridge estimator simultaneously. 
Furthermore, although \cite{nagai2012optimization} proposed 
the model selection criterion of this type, it is not based on estimators 
that are rigorously proven to dominate the maximum likelihood estimator. 

In this paper, we propose new model selection criteria for multivariate linear 
regression based on new shrinkage estimators dominating the maximum likelihood 
estimator under the given risks. 
In particular, even when the model does not include the true model, our proposed estimator dominates the maximum likelihood estimator under the squared risk. 
Moreover, our model selection criteria have the 
following favorite properties: consistency, unbiasedness, and uniformly minimum variance. 
Furthermore, our model selection criteria have closed forms that are given by 
modifying $\aicc$ and $\mcp$.

We construct a class of Bayes estimators that 
dominate the maximum likelihood estimator
under the risks and have a form of the generalized ridge estimator. 
The generalized ridge estimator of multivariate linear regression is a class 
of estimators that can be written as 
\begin{eqnarray*}
\inv{(\ten{\aj}\aj+\Pj\Kj\ten{\Pj})}\ten{\aj}Y,
\end{eqnarray*}
where $\Kj$ is a $\kj \times \kj$ diagonal matrix and $\Pj$ is the $\kj\times\kj$ 
orthogonal matrix of eigenvectors of $\inv{(\ten{\aj}\aj)}$. In other words,
\begin{eqnarray*}
\ten{\Pj}\inv{(\ten{\aj}\aj)}\Pj=\Dj,\m\m \ten{\Pj}\Pj=\I{\kj},
\end{eqnarray*}
where $\Dj=\diag(d_{J,1},d_{J,2},...,d_{J,\kj})$ and  
$d_{J,1}\geq d_{J,2}\geq \cdots\geq d_{J,\kj}$.
Our estimator is related to the generalized Bayes estimator proposed by 
\cite{maruyama2005new} for linear regression, 
the Stein type estimator proposed by \cite{konno1991estimation}
for multivariate linear regression, 
and the generalized Bayes estimator proposed by \cite{tsukuma2009generalized}
for multivariate linear regression.
In contrast to these estimators, our estimators enable us to construct 
closed-form model selection criteria based on them. 
Moreover, as stated above, it is shown that the criteria have 
several favorable statistical properties. 
Since our model selection criteria are based on the 
generalized ridge estimators dominating the maximum likelihood estimator, 
it is expected that the risks of our estimators on the models selected by 
our model selection criteria 
are smaller than the risks of the maximum likelihood estimator on the 
models selected by MCp and AICc. 
We carry out numerical experiments to show the properties of our method.

The contents of this paper are summarized as follows.
In Section 2, we list the classes of estimators dominating 
the maximum likelihood estimator under 
the squared risk and the Kullback-Leibler risk. 
Our estimator is given as a Bayes estimator, and eventually, it is shown that 
it has the same form as a generalized ridge estimator.  
By setting the hyper parameters of our Bayes estimator appropriately, 
we derive the class of estimators dominating the maximum likelihood estimator under the squared risk and the Kullback-Leibler risk.
In Section 3, we construct model selection 
criteria based on the estimators in the classes proposed in Section 2. It is also shown that 
our model selection criteria are uniformly minimum variance unbiased estimators 
of the two risks respectively 
and have consistency. In Section 4, we give numerical 
comparisons between our model selection criteria with AIC, AICc and MCp. 
In Section 5, we provide the discussion and conclusions.

\section{Generalized ridge estimator}
\m\n In this section, we construct a class of Bayes estimators 
that have the form of the generalized ridge estimator and dominate 
the maximum likelihood estimator under the squared risk and 
the Kullback-Leibler risk. 
To derive the estimator, we rotate the coordinate and 
construct the Bayes estimator that can be considered 
as a generalized ridge estimator on the coordinate.
It is shown that the Bayes estimator with tuned hyper parameters dominates 
the maximum likelihood estimator under the squared risk and 
the Kullback-Leibler risk. 
Then, the Bayes estimator is minimax because the maximum likelihood estimator 
is minimax optimal with constant risk.
In particular, in the case of the squared risk, it is not necessary 
that the candidate model includes the true model. 
However, in the case of the Kullback-Leibler risk, 
this property is shown only on models that include the true model.

First, we give a coordinate transformation on $Y$ to derive the estimator.
Let $\qj$ be an $n\times n$ orthogonal matrix such that 
\begin{eqnarray*}
Q_J\aj= \left(
  \begin{array}{ccc}
     \Dj^{-1/2}\ten{\Pj} \\
     0\
  \end{array}
     \right).
\end{eqnarray*}
and let $D_{*J}$ be an $n\times n$ 
diagonal matrix $\diag(d_{J,1},d_{J,2},...,d_{J,\kj},1,...,1)$.
We define random matrices 
$\xj=\ten{(x_{J,1},...,x_{J,\kj})}\in R^{\kj\times p}$ and \\
$\zj=\ten{(z_{J,1},...,z_{J,n-\kj})}\in R^{(n-\kj)\times p}$ such that
\begin{eqnarray*}
\left(
  \begin{array}{ccc}
     \xj \\
     \zj
  \end{array}
     \right)
= D_{*J}^{\frac{1}{2}}\qj Y.
\end{eqnarray*}
Then $(\xj,\zj)$ has the joint density given by
\begin{eqnarray*}
&&(2\pi)^{-\kj p}|\Sigma|^{-\kj}\prod_i^{\kj} \dji^{-p}\exp\left\{-\frac{1}{2}\tr[\Sigma^{-1}\ten{(\xj-\tj)}\Dj^{-1}(\xj-\tj)]\right\}\\
&&\times(2\pi)^{(n-\kj)p}|\Sigma|^{-(n-\kj)}\exp\left\{-\frac{1}{2}\tr[\inv{\Sigma}\sj]\right\}.
\end{eqnarray*}
where $\tj=\ten{\Pj}\bj=\ten{(\theta_{J,1},...,\theta_{J,\kj})}$ and $\sj=\ten{\zj}\zj$. 
With this transformation, the squared risk can be rewritten as 
\begin{eqnarray*}
\E{
\tr\inv{\Sigma}\ten{(\Phi_{J}-\tj)}\Dj(\Phi_{J}-\tj)
},
\end{eqnarray*}
where $\Phi_J$ is an estimator of $\tj$. Note that the 
maximum likelihood estimator of $\tj$ is 
given by $\xj$, and thus $\xj$ is minimax optimal.
We construct the Bayes estimator of $\tj$ 
so that it dominates the maximum likelihood estimator 
under the squared risk and the Kullback-Leibler risk.

\subsection{Derivation of the estimator}
\m\n In this subsection, we derive a Bayes estimator that is based on 
the following prior distribution:
\begin{eqnarray}
\tj|\Sigma&\sim& {\cal N}_{\kj\times p}(0,\Sigma\otimes\Dj(\inv{\Lj}\Cj-\I{\kj}))
=\pi(\tj|\Sigma),\\
\Sigma &\sim& \pi(\Sigma),
\end{eqnarray} 
where $\Lj$ and $\Cj$ are diagonal matrices:  
$\Lj=\diag(\lambda_{J,1},...,\lambda_{J,\kj})$, 
$\Cj=\diag(c_{J,1},...,c_{J,\kj})~(\lji>0,~\cji>0~,i=1,...,\kj)$, 
and $\pi(\Sigma)$ 
is any distribution on the positive definite matrices such that 
$\E{\inv{\Sigma}|\xj,\sj}$ and \\$\inv{\E{\inv{\Sigma}|\xj,\sj}}$ exist. 

\begin{theo}
\n\n
The Bayes estimator based on the prior $(1)$ and $(2)$ is given by
\begin{eqnarray*}
\bjbay=\Pj\tjbay=\Pj(\I{\kj}-\Lj\inv{\Cj})\Dj\ten{\Pj}\ten{\aj}Y
\end{eqnarray*}
and is a generalized ridge estimator.
\end{theo}

{\bf Proof.}\m
To derive the Bayes estimator based on this prior, 
we calculate a part of exponential of the joint density of $(\xj,\zj,\tj,\Sigma)$. 
Let $F_J=\inv{(\inv{\Lj}\Cj-\I{\kj})}$. Then, 
\begin{eqnarray*}
&&\tr\{\Sigma^{-1}\ten{(\xj-\tj)}\Dj^{-1}(\xj-\tj)\} +\tr\{\Sigma^{-1}\ten{\tj}(\Lj^{-1}\Cj-\mathrm{I})^{-1}\Dj^{-1}\tj\}\\
&=&\tr\{\Sigma^{-1}[\ten{(\xj-\tj)}\Dj^{-1}(\xj-\tj)+\ten{\tj}(\Lj^{-1}\Cj-\mathrm{I})^{-1}\Dj^{-1}\tj]\}\\
&=&\tr\{\Sigma^{-1}[\ten{\xj}\inv{\Dj}\xj-\ten{\xj}\inv{\Dj}\tj-\ten{\tj}\inv{\Dj}\xj+
\ten{\tj}(\mathrm{I}+F_J)\inv{\Dj}\tj]\}\\
&=&\tr\left\{\inv{\Sigma}
\left[ 
\ten{\left(\tj-\inv{(\mathrm{I}+F_J)}\xj\right)}(\mathrm{I}+F_J)\inv{\Dj}\left(\tj-\inv{(\mathrm{I}+F_J)}\xj\right)
\right.\right.\\
&&\left.\left.~~~~~~~~~~~~~~~~~~~~~~+\ten{\xj}\inv{\Dj}\left(\mathrm{I}-\inv{(\mathrm{I}+F_J)}\right)\xj
\right]\right\}, 
\end{eqnarray*}
and $\inv{(\I{\kj}+F_J)}$ is given by
\begin{eqnarray*}
\inv{(\I{\kj}+F_J)}&=&\inv{\left(\I{\kj}+\inv{\left(\Cj\inv{\Lj}-\I{kj}\right)}\right)}\\
&=&\inv{\left(\I{\kj}+\Lj\inv{\left(\Cj-\Lj\right)}\right)}\\
&=&\inv{(\Cj\inv{(\Cj-\Lj)})}\\
&=&\I{\kj}-\Lj\inv{\Cj}.
\end{eqnarray*}
Therefore, the joint density of $(\xj,\zj,\tj,\Sigma)$ is proportional to
\begin{eqnarray*}
&&|\Sigma|^{-\frac{n}{2}-\frac{p}{2}}
\mathrm{exp}\left\{
-\frac{1}{2}\tr\left\{\inv{\Sigma}
\left[ 
\ten{\left(\tj-(I-\Lj\inv{\Cj})\xj\right)}\inv{( I-\Lj\inv{\Cj})}\inv{\Dj}
\right. \right. \right. \\
&&\m\m\m\m\times \left.  \left. \left. 
\left(\tj-( I-\Lj\inv{\Cj})\xj\right) +
\ten{\xj}\inv{\Dj}\left(I- (I-\Lj\inv{\Cj})\right)\xj
\right]\right\} \right.\\
&&\m\m\m\m
-\left.\frac{1}{2}\tr[\inv{\Sigma}\sj]
\right\}\pi(\Sigma).
\end{eqnarray*}
The Bayes estimator of $\Theta$ under the squared risk is obtained by 
\begin{eqnarray*}
&&\frac{\partial}{\partial \Phi_J}\E
{\tr\{\inv{\Sigma}\ten{(\tj-\Phi_J)}\Dj^{-a}(\tj-\Phi_J)\}|\xj,\sj
}=0\\
\Rightarrow&&\E
{2\Dj^{-a}(\tj-\Phi_J)\inv{\Sigma}|\xj,\sj
}=0\\
\Rightarrow&&\E{\Dj^{-a}\tj\inv{\Sigma}|\xj,\sj}=\E{\Dj^{-a}\Phi_J\inv{\Sigma}|\xj,\sj}\\
\Rightarrow&&\E{\tj\inv{\Sigma}|\xj,\sj}=\Phi_J\E{\inv{\Sigma}|\xj,\sj}\\
\Rightarrow&&\Phi_J=\E{\tj\inv{\Sigma}|\xj,\sj}\inv{\{\E{\inv{\Sigma}|\xj,\sj}\}}\\
\Rightarrow&&\Phi_J=\E{\tj|\xj,\sj}.
\end{eqnarray*}
Therefore, the Bayes estimator $\tjbay$ based on this prior is $\E{\tj|\xj,\sj}$. 
From the joint density of $(\xj,\zj,\tj,\Sigma)$, 
\begin{eqnarray*}
\tjbay=(\I{\kj}-\Lj\inv{\Cj})\xj.
\end{eqnarray*} 
Moreover, the Bayes estimator $\bjbay$ of $\bj$ based on this prior can 
be written as 
\begin{eqnarray*}
\bjbay&=&\Pj\tjbay=\Pj(\I{\kj}-\Lj\inv{\Cj})\Dj\ten{\Pj}\ten{\aj}Y\\
&=&\Pj\inv{(\inv{\Dj}\inv{(\I{\kj}-\Lj\inv{\Cj})})}\ten{\Pj}\ten{\aj}Y\\
&=&\Pj\inv{(\inv{\Dj}(\I{\kj}+\Lj\inv{(\Cj-\Lj)})}\ten{\Pj}\ten{\aj}Y\\
&=&\Pj\inv{(\inv{\Dj}+\inv{\Dj}\Lj\inv{(\Cj-\Lj)})}\ten{\Pj}\ten{\aj}Y.
\end{eqnarray*}
Furthermore, let $K_J=\inv{\Dj}\Lj\inv{(\Cj-\Lj)}$. Then $\bjbay$ is regarded as the 
generalized ridge estimator. \qed

The Bayes estimator is close to the multivariate form of Maruyama (2005). 
While he defined a prior of $\Lj$ as univariate, in this paper,
$\Lj$ is a hyper parameter.
Since his estimator does not have a closed form, 
we change the prior to construct 
the Bayes estimator expressed in a closed form. 
We consider a plug-in predictive density by using this 
Bayes estimator even for the Kullback-Leibler risk because it is easy to handle.

\subsection{Generalized ridge estimator under the squared risk}
\m\n In this subsection, we set the parameters of the generalized ridge estimator 
and show that it dominates the maximum likelihood estimator under 
the squared risk for any candidate 
model. In many studies, although several generalized ridge estimators have been examined, 
they assume that 
the candidate model includes the true model. 
However, this paper does not make this assumption.

In the following theorem, we give sufficient conditions of the parameters 
so that the Bayes estimator 
$\bjbay$ dominates the maximum likelihood estimator $\bjmle$ under the squared risk.

\begin{theo}
\n\n
\label{thmSq1}
Let
$\Cj=\diag(\ten{x_{J,1}}\inv{S_F}x_{J,1},...,
\ten{x_{J,\kj}}\inv{S_F}x_{J,\kj})$ 
and assume that
 \begin{description} 
\item[$(i)$] $0 < \lji < \frac{2\dji(p-2)}{n-p-\kf+3}~(i=1,...,\kj),$
\item[$(ii)$] $p\geq3$ and $n-p-\kf+3>0,$
\item[$(iii)$] $J_*\in {\cal J}.$
 \end{description} 
Then, 
\begin{eqnarray*}
\srisk{\bjt}{\Sigma}{\hat{B}_{J\mathrm{B}}} < \srisk{\bjt}{\Sigma}{\hat{B}_J}\m 
({}^{\forall}J\in{\cal J}).
\end{eqnarray*}
Moreover $\srisk{\bjt}{\Sigma}{\hat{B}_{J\mathrm{B}}}$ is minimum at 
$\lji=\frac{\dji(p-2)}{n-p-\kf+3}~(i=1,...,\kj)$.
\end{theo}

We employed $\cji=\ten{\xjic}\inv{S_F}\xjic~(i=1,2,...,\kj)$
instead of $\ten{\xjic}\inv{S_J}\xjic$ 
because of the following reason. 
If we use $S_J$ and not $S_F$, we need to assume that 
the candidate model $J$ includes the true model to 
show that $\hat{B}_{J,B}$ dominates $\hat{B}_J$. 
Moreover, in this case, we do not obtain 
a closed form of the model selection criterion. 
However, by employing our definition of $\Cj$, 
we obtain an unbiased estimator of the risk by slightly modifying $\mcp$.

While the assumption $(i)$ of this theorem gives the condition of $\Lj$,
we may simply set $\lji = \dji(p-2)/(n-\kj-p+3)~(i=1,2,...,\kj)$. 
The assumption $(ii)$ of this theorem means a restriction of the dimension. 
In particular, $p\geq3$ is the same condition 
as that under which Stein's estimator dominates 
the maximum likelihood estimator with respect to the squared risk. 
The assumption $(iii)$ is 
to ensure that the candidate models contain the true model.

Each row of $\bjbay$ is similar to the Stein's estimator. Furthermore, $\bjbay$ 
dominates $\bjmle$ under the squared risk when 
the model $J$ includes the true model.
However, this is not obvious when the model $J$ does not include the true model. 
Furthermore, from the form of lows of $\bjbay$ and the fact that 
Stein's estimator is not admissible, $\bjbay$ is not admissible.
However, we employ this form rather than pursuing a statistically optimal one 
because of its computational efficiency. In fact, the model selection criteria 
based on the estimator can be analytically computed as shown later.\\

{\sc Proof.}\m
From the assumption $(iii)$, without loss of generality, we may regard 
$J_*$ as $F$.
Therefore, let $\bjt$ be $B_F$ and $\ajt$ be $\af$. For ${}^{\forall}J\in{\cal J}$,
\begin{eqnarray*}
\srisk{\bjt}{\Sigma}{\bjbay}&=&\E{
\tr\left[\inv{\Sigma}\ten{(\bjbay-\bjt)}(\ten{\af}\af)(\bjbay-\bjt)\right]
}\\
&=&\E{
\tr\left[\inv{\Sigma}\ten{(\aj\bjbay-\af\bjt)}(\aj\bjbay-\af\bjt)\right]
}\\
&=&\srisk{\bjt}{\Sigma}{\bjmle}\\
&&-2\E{
\tr\left[\inv{\Sigma}\ten{(\aj\bjmle-\af\bjt)}
\aj\Pj\Lambda_J\Cj^{-1}\ten{\Pj}\bjmle\right]
}\\
&&+\E{\tr\left[\inv{\Sigma}\ten{\bjmle}\Pj\Lambda_J\Cj^{-1}\ten{\Pj}\ten{\aj}
\aj\Pj\Lambda_J\Cj^{-1}\ten{\Pj}\bjmle\right]
}.
\end{eqnarray*}
Therefore it is sufficient to show that the sum of the second and third terms 
is non-positive. From the definitions of $\xj$ and $S_F$,
\begin{eqnarray*}
\xj\sim {\cal N}_{\kj\times p}(\Dj\ten{\Pj}\ten{\aj}\af\bjt,\Sigma\otimes \Dj),~~~~
S_F \sim W_p(n-\kf,\I{\kj}).
\end{eqnarray*}
Let $\wj=\Dj^{-\frac{1}{2}}\xj=\ten{(w_{J,1},...,w_{J,\kj})}$. Then 
\begin{eqnarray*}
\wj&=&\ten{\Pj}(\ten{\aj}\aj)^{\frac{1}{2}}\bjmle\sim {\cal N}_{\kj\times p}(\Psi_J,\Sigma\otimes I),\\
\ten{\xjic}\isf\xjic&=&\dji\ten{\wjic}\isf\wjic,
\end{eqnarray*}
where $\Psi_J=\ten{\Pj}(\ten{\aj}\aj)^{-\frac{1}{2}}\aj\af\bjt=\ten{(\psi_{J,1},
...,\psi_{J,\kj})}$.
From the definition of $\wj$ and the joint density of $(\xj,\sj)$, the joint density 
of \\$w_{J,1},w_{J,2},...,w_{J,\kj},\zf$, and
$\ten{Y}(\af(\ten{\af}\af)^{-1}\ten{\af}-\aj(\ten{\aj}\aj)^{-1}\ten{\aj})Y)$ is given by 
\begin{eqnarray*}
&&(2\pi)^{-\kj p}|\Sigma|^{-\kj}\prod_{i=1}^{\kj}\exp\left\{-\frac{1}{2}
\ten{(\wjic-\psi_{J,i})}\inv{\Sigma}(\wjic-\psi_{J,i})\right\}\\
&&\times(2\pi)^{(n-\kj)p}|\Sigma|^{-(n-\kj)}\exp\left\{
-\frac{1}{2}\tr[\inv{\Sigma}S_F] \right. \\
&&\m\m\m\m\m\m\m- \left. \frac{1}{2}\tr[\inv{\Sigma}
\ten{Y}(\af(\ten{\af}\af)^{-1}\ten{\af}-\aj(\ten{\aj}\aj)^{-1}\ten{\aj})Y]\right\}.
\end{eqnarray*}
From the Fisher-Cochran theorem and the definitions of $\wj$ and  $S_F$,  $\wjic~(i=1,2,...,\kj)$, $S_F$, and $\ten{Y}(\af(\ten{\af}\af)^{-1}\ten{\af}-\aj(\ten{\aj}\aj)^{-1}\ten{\aj})Y$ 
are independent for a fixed model $J$. 
Therefore, from Lemma 2.1 of \cite{kubokawa2001robust}, the second term is given by
\begin{eqnarray*}
&&-2\E{
\tr\left[\inv{\Sigma}\ten{\left(\aj\bjmle-\af\bjt\right)}
\aj\Pj\Lambda_J\Cj^{-1}\ten{\Pj}\bjmle\right]
}\\
&=&-2\E{
\tr\left[\inv{\Sigma}\ten{\left(\ten{\aj}\aj\bjmle-\ten{\aj}\af\bjt\right)}
\Pj\Lambda_J\Cj^{-1}\ten{\Pj}\bjmle\right]
}\\
&=&-2\E{
\tr\left[\inv{\Sigma}\ten{\left(\left(\ten{\aj}\aj\right)^{\frac{1}{2}}\bjmle-\left(\ten{\aj}\aj\right)^{-\frac{1}{2}}\aj\af\bjt\right)}
\left(\ten{\aj}\aj\right)^{\frac{1}{2}}\right.\right.\\
&&~~~~~~~~\times\left.\left.\Pj\Lambda_J\Cj^{-1}\ten{\Pj}\left(\ten{\aj}\aj\right)^{-\frac{1}{2}}\left(\ten{\aj}\aj\right)^{\frac{1}{2}}\bjmle\right]
}\\
&=&-2\E{
\tr\left[\inv{\Sigma}\ten{\left(W_J-\ten{\Pj}\left(\ten{\aj}\aj\right)^{-\frac{1}{2}}\aj\af\bjt\right)}
\Dj^{-\frac{1}{2}}\Lambda_J\Cj^{-1}\Dj^{\frac{1}{2}}W_J\right]
}\\
&=&-2\E{
\tr\left[\inv{\Sigma}\ten{\left(\wj-\ten{\Pj}\left(\ten{\aj}\aj\right)^{-\frac{1}{2}}\aj\af\bjt\right)}
\Lambda_J\Cj^{-1}\wj\right]
}\\
&=&-2\sum_{i=1}^{\kj}\E{
\lji\inv{\cji}\ten{\left(W_J-\ten{\Pj}\left(\ten{\aj}\aj\right)^{-\frac{1}{2}}\aj\af\bjt\right)}_{i\cdot}\inv{\Sigma}\wji
}\\
&=&-2\sum_{i=1}^{\kj}\lji\inv{\dji}\\
&&\m\times\E{
\invwjicsfwjic\ten{\left(\wj-\ten{\Pj}\left(\ten{\aj}\aj\right)^{-\frac{1}{2}}\aj\af\bjt\right)}_{i\cdot}\inv{\Sigma}\wji
}\\
&=&-2\sum_{i=1}^{\kj}\lji\inv{\dji}\E{
\tr\left(\nabla_{i\cdot}\invwjicsfwjic\wjic\right)
}\\
&=&-2\sum_{i=1}^{\kj}\sum_{j=1}^{p}\lji\inv{\dji}\E{
\frac{\partial}{\partial W_{J,ij}}\invwjicsfwjic W_{J,ij}
}\\
&=&-2p\sum_{i=1}^{\kj}\lji\inv{\dji}\E{
\invwjicsfwjic
}\\
&&+4\sum_{i=1}^{\kj}\sum_{j=1}^{p}\lji\inv{\dji}\E{
(\wjicsfwjic)^{-2}
\left(\ten{\wjic}\isf\right)_jW_{J,ij}
}\\
&=&-2p\sum_{i=1}^{\kj}\E{
\lji\inv{\cji}
}+4\sum_{i=1}^{\kj}\E{\lji\cji^{-1}
}\\
&=&-2\left(p-2\right)\sum_{i=1}^{\kj}\E{
\lji\inv{\cji}
}.
\end{eqnarray*}
Similarly, from the proof of Proposition 2.1 of \cite{kubokawa2001robust}, the third term is given by
\begin{eqnarray*}
&&\E{
\tr\left(\inv{\Sigma}\ten{\bjmle}\Pj\Lambda_J\Cj^{-1}\ten{\Pj}\ten{\aj}
\aj\Pj\Lambda_J\Cj^{-1}\ten{\Pj}\bjmle\right)
}\\
&=&\E{
\tr\left(\inv{\Sigma}\ten{\xj}\Lambda_J\Cj^{-1}\Dj^{-1}\Lambda_J\Cj^{-1}\xj\right)
}\\
&=&\E{
\tr\left(\inv{\Sigma}\ten{\wj}\Lj^2\Cj^{-2}\wj\right)
}\\
&=&\sum_{i=1}^{\kj}\lji^2\dji^{-2}\E{
(\wjicsfwjic)^{-2}\ten{\wjic}\inv{\Sigma}\wjic
}\\
&=&\sum_{i=1}^{\kj}\lji^2\dji^{-2}\E{
\invwjicsfwjic
}\\
&=&(n-\kf-p+3)\sum_{i=1}^{\kj}\E{
\inv{\dji}\lji^2\inv{\cji}
}.
\end{eqnarray*}
Therefore, from  the assumptions $(i)$ and $(ii)$, the sum of the second and third 
terms is 
given by
\begin{eqnarray*}
&&-2\left(T-2\right)\sum_{i=1}^{\kj}\E{
\lji^2\inv{\cji}
}+(n-\kf-T+3)\sum_{i=1}^{\kj}\E{
\inv{\dji}\lji^2\inv{\cji}
}\\
&=&\sum_{i=1}^{\kj}\lji\E{\inv{\cji}}\left\{-2(T-2)+(n-\kf-T+3)\lji\dji^{-1}\right\}\\
&<&0.
\end{eqnarray*}
Therefore, $\srisk{\bjt}{\Sigma}{\bjbay}$ is smaller than $\srisk{\bjt}{\Sigma}{\bjmle}$. 
Moreover, \\$\srisk{\bjt}{\Sigma}{\bjbay}$ is minimum at 
$\lji = \frac{\dji(p-2)}{(n-\kj-p+3)}~(i=1,2,...,\kj)$ because $\srisk{\bjt}{\Sigma}{\bjmle}$ is constant.
\qed\\

\subsection{Generalized ridge estimator under the Kullback-Leibler risk}
\m\n In this subsection, we set the parameters of the generalized ridge estimator 
and show that it dominates the maximum likelihood estimator 
under the Kullback-Leibler risk
when the candidate model includes the true model. We consider 
a plug-in predictive density that is obtained by plugging-in estimators
to $B_J$ and $\Sigma$
to construct a model selection criterion 
which is given in a closed form. 

Let $\hat{\Sigma}_J$ be the maximum likelihood 
estimator of $\Sigma$ on the model $J$. Then 
we obtain the condition of the parameters under which the plug-in predictive density 
with $\bjbay$ and $\hat{\Sigma}_J$ dominates the plug-in predictive density 
with $\bjmle$ and $\hat{\Sigma}_J$ under the Kulback-Leibler risk. 

\begin{theo}
\n\n
\label{thmKL}
Let $\Cj=\diag(\ten{x_{J,1\cdot}}\inv{S_{J}}x_{J,1\cdot},...,\ten{x_{J,\kj\cdot}}\inv{S_{J}}x_{J,\kj\cdot})$ and assume that
 \begin{description} 
\item[$(i)$] $0<\lji<\frac{2\dji(p-2)}{n-\kj-p+1},$
\item[$(ii)$] $p \geq 3$ and $n-p-\kj-1>0,$
\item[$(iii)$] $J_* \subset J.$
 \end{description} 
Then,
\begin{eqnarray*}
\klrisk{\bjt}{\Sigma}{f(\cdot|\bjbay,\hat{\Sigma}_J)}<
\klrisk{\bjt}{\Sigma}{f(\cdot|\bjmle,\hat{\Sigma}_J)}.
\end{eqnarray*}
Moreover, $R_{KL}(\bjt,\Sigma,f(\cdot|\bjbay,\hat{\Sigma}_J))$ is minimum 
at $\lji=\frac{\dji(p-2)}{n-\kj-T+1}~(i=1,...,\kj)$.
\end{theo}

In the case of the Kullback-Leibler risk, we exclude the case where 
the candidate model does not include the true model. 
The reason is as follows. 
Our estimator is not of the covariance but the mean. However, the plug-in predictive 
density depends on not only the mean but also the covariance and thus the predictive 
performance is affected by the estimation performance of covariance.
This makes it difficult to analyze whether the plug-in predictive density
dominates the maximum likelihood estimator in the case where $J_* \not\subset J$.
Even for this situation, it might be possible to construct a plug-in predictive density 
that dominates the plug-in predictive density with the maximum likelihood estimator.
However, our main purpose is to construct a model selection criterion, 
and thus, we do not pursue this problem in this paper. \\

{\sc Proof.}\m From the assumption $(iii)$,   
without loss of generality, we regard $J_*$ as $J$.
The Kullback-Leibler risk of the plug-in predictive 
density based on $\bjbay$ and $\Sigmajmle$ is given by
\begin{eqnarray*}
\klrisk{\bj}{\Sigma}{f(\cdot|\bjbay,\Sigmajmle)}=
\Es{\tilde{Y},Y}{\log{f(\tilde{Y}|\bj,\Sigma)}}
-\Es{\tilde{Y},Y}{\log{f(\tilde{Y}|\bjbay,\hat{\Sigma}_J)}}.
\end{eqnarray*}
Thus, the first terms depends on only the true distribution and not on the plug-in
predictive density.
The integrand of the second term can be written as
\begin{eqnarray*}
-\frac{n}{2}\log{|\hat{\Sigma}_J|}-\frac{np}{2}\log{2\pi}
-\frac{1}{2}\tr\left\{\inv{\hat{\Sigma}_J}\ten{(\tilde{Y}-\aj\bjbay)}(\tilde{Y}-\aj\bjbay)\right\}.
\end{eqnarray*}
The third term of this expression can be written as
\begin{eqnarray*}
&&\tr\left\{\inv{\hat{\Sigma}_J}\ten{(\tilde{Y}-\aj\bjmle)}(\tilde{Y}-\aj\bjmle)\right\}\\
&+&2\tr\left\{\inv{\hat{\Sigma}_J}\ten{(\tilde{Y}-\aj\bjmle)}(\aj\bjmle-\aj\bjbay)\right\}\\
&+&\tr\left\{\inv{\hat{\Sigma}_J}\ten{(\aj\bjbay-\aj\bjmle)}(\aj\bjbay-\aj\bjmle)\right\}.
\end{eqnarray*}
Similarly, $\log{f(\tilde{Y}|\bjmle,\hat{\Sigma}_J)}$ can be written as
\begin{eqnarray*}
-\frac{n}{2}\log{|\hat{\Sigma}_J|}-\frac{np}{2}\log{2\pi}
-\frac{1}{2}\tr\left\{\inv{\hat{\Sigma}_J}\ten{(\tilde{Y}-\aj\bjmle)}(\tilde{Y}-\aj\bjmle)\right\}.
\end{eqnarray*}
Therefore,
\begin{eqnarray*}
&&2R_{\mathrm{KL}}(\bj,\Sigma,f(\cdot|\bjbay,\hat{\Sigma}_J))
-2R_{\mathrm{KL}}(\bj,\Sigma,f(\cdot|\bjmle,\hat{\Sigma}_J))\\
&=&2\Es{\tilde{Y},Y}{
\log{f(\tilde{Y}|\bjmle,\hat{\Sigma}_J)}-\log{f(\tilde{Y}|\bjbay,\hat{\Sigma}_J)}
}\\
&=&2\Es{\tilde{Y},Y}{
\tr\left(\inv{\hat{\Sigma}_J}\ten{(\tilde{Y}-\aj\bjmle)}(\aj\bjmle-\aj\bjbay)\right)
}\\
&&\m+\Es{\tilde{Y},Y}{\tr\left(\inv{\hat{\Sigma}_J}
\ten{(\aj\bjbay-\aj\bjmle)}(\aj\bjbay-\aj\bjmle)\right)
}.
\end{eqnarray*}
The second term of the last expression above is evaluated by 
\begin{eqnarray*}
&&\tr\left\{\inv{\hat{\Sigma}}\ten{(\aj\bjbay-\aj\bj)}(\aj\bjbay-\aj\bj)\right\}\\
&=&n\E{\tr(\isj\ten{\xj}\Lj^2\Cj^{-2}\Dj^{-1}\xj)}\\
&=&n\sum_{i=1}^{\kj}\lji^2\dji^{-1}\E{\invjxisx}.
\end{eqnarray*}
Let $\wj=\Dj^{-\frac{1}{2}}\xj=\ten{(w_{J,1},...,w_{J,\kj})}\sim
{\cal N}_{\kj\times p}(\Dj^{-\frac{1}{2}}\tj,\Sigma\otimes\I{\kj})$.
Then from Lemma 2.1 of \cite{kubokawa2001robust}, 
the first term is given by
\begin{eqnarray*}
&&2\Es{\tilde{Y},Y}{
\tr\left(\inv{\hat{\Sigma}_J}\ten{(\tilde{Y}-\aj\bjmle)}(\aj\bjmle-\aj\bjbay)\right)
}\\
&=&2\Es{Y}{
\tr\left(\inv{\hat{\Sigma}_J}\ten{(\aj\bj-\aj\bjmle)}(\aj\bjmle-\aj\bjbay)\right)
}\\
&=&2\Es{Y}{
n\tr\left(\isj\ten{(\Dj^{-\frac{1}{2}}\tj-\wj)}\Lj \Cj^{-1}\wj\right)
}\\
&=&-2\Es{Y}{
n\sum_{i=1}^{\kj}\lji\cji^{-1}\ten{(\wjic-\dji^{-\frac{1}{2}}\tjic)}\isj\wjic
}\\
&=&-2n\sum_{i=1}^{\kj}\lji\tr\left\{\Es{Y}{
\cji^{-1}\isj\wjic\ten{(\wjic-\dji^{-\frac{1}{2}}\tjic)}
}\right\}\\
&=&-2n\sum_{i=1}^{\kj}\lji\inv{\dji}\Es{Y}{
\tr\left(\nabla_i \left(\invwjicsjwjic\isj\wjic\right)\Sigma\right)
}\\
&=&-2n\sum_{i=1}^{\kj}\lji\inv{\dji}\Es{Y}{
\sum_{j=1}^p\sum_{k=1}^p
\frac{\partial}{\partial W_{J,ij}}\left(\invwjicsjwjic(\isj\wjic)_k\right)\Sigma_{jk}
}\\
&=&-2n\sum_{i=1}^{\kj}\lji\inv{\dji}\Es{Y}{
\sum_{j=1}^p\sum_{k=1}^p\left(
-2(\wjicsjwjic)^{-2}(\wjic\isj)_j\Sigma_{jk}(\isj\wjic)_k\right.\right.\\
&&\m\m\m\m\m\m\m\m\m\m\m\m\m+\left.\left.
\invwjicsjwjic(\inv{\sj})_{kj}\Sigma_{jk}
\right)}\\
&=&-2n\sum_{i=1}^{\kj}\lji\inv{\dji}\Es{Y}{
-2(\wjicsjwjic)^{-2}\ten{\wjic}\isj\Sigma\isj\wjic \right. \\
&&\left.\m\m\m\m\m\m\m\m\m\m\m\m\m+\invwjicsjwjic\tr(\isj\Sigma)
}.
\end{eqnarray*}
Let $\wjic'=\Sigma^{-\half}\wjic$, $\zj'=\zj\Sigma^{-\half}$ and 
$\sj'=\ten{\zj'}\zj'$. Then\\ 
$\wjic'\sim {\cal N}_T(\dji^{-\frac{1}{2}}\Sigma^{-\frac{1}{2}}\tjic,\mathrm{I})$,
$\zj' \sim {\cal N}_{(n-\kj)\times p}$, and 
$\sj'\sim {\cal W}_p(n-\kj,\mathrm{I})$ where 
${\cal W}_p(k,\Sigma)$ is the Wishart distribution that has 
degree of freedom $k$ and scale matrix $\Sigma$. 
Let $\Rji{1}$ be a $p\times p$ orthogonal matrix such that \\
$\Rji{1}\wjic'=\ten{(\sqrt{\ten{\wjic'}\wjic'},0,...,0)}$, and let 
$\Rji{1}\ten{\zj}=\ten{(\vji{1},\vji{2})}$ where $\vji{1}$ is an $(n-\kj)\times 1$ vector 
and $\vji{2}$ is an $(n-\kj)\times(p-1)$ matrix. 
Furthermore, let $\Rji{2}$ be an $(n-\kj)\times(n-\kj)$ orthogonal matrix such 
that $\Rji{2}\vji{2}=(0,(\ten{\vji{2}}\vji{2})^{\frac{1}{2}})$, and let 
$\Rji{2}\vji{1}=\ten{(\ten{\uji{1}},\ten{\uji{2}})}$ where $\uji{1}$ is a $(n-p-\kj+1)\times 1$ 
vector and $\uji{2}$ is a $(p-1)\times 1$ vector. Then
\begin{eqnarray*}
\wjic'\sj'^{-1}\wjic'&=&\frac{\ten{\wjic'}{\wjic'}}{\ten{\vji{1}}(\I{n-\kj}
-\vji{2}(\ten{\vji{2}}\vji{2})^{-1}\ten{\vji{2}})\vji{1}}\\
&=&\frac{\ten{\wjic'}{\wjic'}}{\ten{\uji{1}}\uji{1}}
\end{eqnarray*}
and let 
\begin{eqnarray*}
V_{J,i}&=&\Rji{1}\ten{\zj}\ten{\zj}\Rji{1}
=\Rji{1}\sj\Rji{1}\\
&=&\left(
\begin{array}{ccc}
\Vji{11}{} & \Vji{12}{}\\
\Vji{21}{} & \Vji{22}{}
\end{array}
\right)
=\left(
\begin{array}{ccc}
\vjivji{1}{1} & \vjivji{1}{2}\\
\vjivji{2}{1} & \vjivji{2}{2}
\end{array}
\right).
\end{eqnarray*}
Then 
\begin{eqnarray*}
&&\wjic'\sj'^{-2}\wjic'\\
&=&\ten{\wjic'}\wjic'(1,0,...,0)\Vji{}{-2}\ten{(1,0,...,0)}\\
&=&\ten{\wjic'}\wjic'((\Vji{}{-1})_{11}^2+(\Vji{}{-1})_{12}(\Vji{21}{-1}))\\
&=&\ten{\wjic'}\wjic'(\Vji{11\cdot2}{-2}+\Vji{11\cdot2}{-2}\Vji{12}{}\Vji{22}{-2}
\Vji{21}{})\\
&=&\frac{\ten{\wjic'}\wjic'}{(\ten{\uji{1}}\uji{1})^2}
(1+\vjivji{1}{2}(\vjivji{2}{2})^{-2}\vjivji{2}{1})\\
&=&\frac{\ten{\wjic'}\wjic'}{(\ten{\uji{1}}\uji{1})^2}
(1+\ten{\uji{2}}(\vjivji{2}{2})^{-1}\uji{2}),
\end{eqnarray*}
where $\Vji{11\cdot2}{}=\Vji{11}{}-\Vji{12}{}\Vji{22}{-1}\Vji{21}{}$. 
From the definitions of $\wjic$, $\uji{1}$, $\uji{2}$ and $\vji{2}$, they are 
independent. Therefore, 
\begin{eqnarray*}
&&-2n\sum_{i=1}^{\kj}\lji\inv{\dji}\E{
-2(\wjicsjwjic)^{-2}\ten{\wjic}\isj\Sigma\isj\wjic\right.\\
&&\left.\m\m\m\m\m\m\m\m\m+\invwjicsjwjic\tr(\isj\Sigma)
}\\
&=&-2n\sum_{i=1}^{\kj}\lji\inv{\dji}\E{
-2(\wjic'\sj'^{-1}\wjic')^{-2}\ten{\wjic'}\sj'^{-2}\wjic'
+\inv{\wjic'\sj'^{-1}\wjic'}\tr\sj'^{-1}
}\\
&=&-2n\sum_{i=1}^{\kj}\lji\inv{\dji}\E{
-2\frac{1}{\ten{\wjic'}\wjic'}(1+\ten{\uji{2}}(\vjivji{2}{2})^{-1}\uji{2})\right.\\
&&\m\m\m\m\m\m\m\m\m\m\m
+\left.\frac{\ten{\uji{1}}\uji{1}}{\ten{\wjic'}\wjic'}
((\ten{\uji{1}}\uji{1})^{-1}+\tr(\Vji{}{-1})_{22})
}\\
&=&-2n\sum_{i=1}^{\kj}\lji\inv{\dji}\E{
\frac{1}{\ten{\wjic'}\wjic'}
}\\
&&\m\times
\left\{-1+\E{
\ten{\uji{1}}\uji{1}
}\E{\tr(\Vji{}{-1})_{22}
}-2\E{\ten{\uji{2}}(\vjivji{2}{2})^{-1}\uji{2}
}\right\}\\
&=&-2n\sum_{i=1}^{\kj}\frac{\lji\inv{\dji}}{n-\kj-p+1}\E{
\frac{\ten{\uji{1}}\uji{1}}{\ten{\wjic'}\wjic'}
}\\
&&\m\m\m\m\m\m\times
\left\{-1+\frac{(n-\kj-p+1)(p-1)}{n-\kj-p-1}-\frac{2(p-1)}{n-\kj-p-1}
\right\}\\
&=&-2n\sum_{i=1}^{\kj}\frac{\lji\inv{\dji}(p-2)}{n-\kj-p+1}\E{
\invwjicsjwjic
}\\
&=&-2n\sum_{i=1}^{\kj}\frac{\lji(p-2)}{n-\kj-p+1}\E{
\inv{\cji}.
}
\end{eqnarray*}
Therefore from the assumptions $(i)$ and $(ii)$, this implies that
\begin{eqnarray*}
&&2\klrisk{\bj}{\Sigma}{f(\cdot|\bjbay,\hat{\Sigma}_J)}
-2\klrisk{\bj}{\Sigma}{f(\cdot|\bjmle,\hat{\Sigma}_J)}\\
&=&n\sum_{i=1}^{\kj}\lji\E{\cji^{-1}}\left(
-\frac{2(p-2)}{n-\kj-p+1}+\dji^{-1}\lji\right)\\
&<&0.
\end{eqnarray*}
Furthermore, it is obvious that $\klrisk{\bj}{\Sigma}{f(\cdot|\bjbay,\hat{\Sigma}_J)}$
is minimum at\\
$\lji=\frac{\dji(p-2)}{n-\kj-p+1}~(i=1,...,\kj)$.
\qed\\

\section{Model selection criterion}
\m\n In this section, we construct model selection criteria based on the generalized ridge 
estimators that are proposed in Section 2. We show that 
the model selection criteria are 
unbiased estimators of the risks of the generalized ridge estimators.
Moreover, we show that they are
uniformly minimum variance unbiased 
estimators and have consistency.

We consider a noncentrality matrix to show consistency. Let $\rj=\kf-\kj$ and
\begin{eqnarray*}
\tilde{\Omega}_J=\Sigma^{-\half}\ten{(\ajt\bjt)}(\af(\ten\af\af)^{-1}\ten\af
-\aj(\ten\aj\aj)^{-1}\ten\aj)\ajt\bjt\Sigma^{-\half}.
\end{eqnarray*}
Then, we can express $\tilde{\Omega}_J=\ten\Gamma_J\Gamma_J$ where
$\Gamma_J$ is an $\rj\times p$ matrix because the rank of 
$\tilde{\Omega}_J$ is at most $\rj$. Moreover, let
\begin{eqnarray*}
\Omega_J&:=&\Gamma_J\ten\Gamma_J,\\
\Xi_J&:=&\frac{1}{np}\Omega_J,
\end{eqnarray*}
where $\Omega_J$ and $\Xi_J$ are $\rj \times\rj$ matrices. We call $\Omega_J$ 
the noncentrality matrix of $Y$ on $J$ and 
let the rank of $\Omega_J$ be denoted by $\gamma_J$. We assume that 
$\gamma_J$ is independent of $n$ and $p$.
Intuitively, $\Omega_J$ represents ``magnitude'' of model misspecification. 
Indeed, it holds that 
\begin{eqnarray}
\Omega_J= 0~~~(\forall J\supset J_*)
\end{eqnarray}
because $\af\inv{(\ten\af\af)}\ten\af$ and $\aj\inv{(\ten\aj\aj)}\ten\aj$ are 
projection matrices where the range of 
$\af\inv{(\ten\af\af)}\ten\af-\aj\inv{(\ten\aj\aj)}\ten\aj$
is perpendicular to the range of $\ajt$. 

\subsection{Modified MCp}
\m\n In this subsection, we propose a model selection criterion 
that is a modification of $\mcp$ under the squared risk. 
$\mcp$ is the estimator of the squared risk with the maximum likelihood estimator 
and given by
\begin{eqnarray*}
\mcp(J)=(n-\kf-p-1)\tr(\inv{S_F}S_J)+p(2\kj+p+1-n),
\end{eqnarray*}
where $\mcp(J)$ is $\mcp$ under a model $J$.
We construct a model selection criterion based on the 
generalized ridge estimator dominating the maximum likelihood estimator.
Let $\Lj=\diag\left(\frac{d_{J,1}(p-2)}{n-\kf-p+3},...,\frac{d_{J,\kj}(p-2)}{n-\kf-p+3}\right)$ 
and $\Cj = \diag\left(\ten{x_{J,1}}\inv{S_F}x_{J,1},...,
\ten{x_{J,\kj}}\inv{S_F}x_{J,\kj}\right)$.
The risk is minimum with this setting. 
Then
\begin{eqnarray*}
\srisk{\bjt}{\Sigma}{\bjbay}=\srisk{\bjt}{\Sigma}{\bjmle}-(p-2)\E{
\tr(\Lj\inv{\Cj})}
\end{eqnarray*}
by the proof of Theorem 2.
Based on this observation, we propose to use $\zmcp$ as an unbiased estimator of 
$\srisk{\bjt}{\Sigma}{\bjbay}$ under a model $J$, which is defined as
\begin{eqnarray*}
&&\zmcp(J) \\
&=& \mcp(J) -(p-2)\tr(\Lj\inv{\Cj})\\
&=&(n-\kf-p-1)\tr(\inv{S_F}S_J)+p(2\kj+p+1-n)
-(p-2)\tr(\Lj\inv{\Cj}).
\end{eqnarray*}

The model selection criterion has the following properties.
\begin{theo}
\n\n
$\zmcp$ is a uniformly minimum variance unbiased estimator 
of the squared risk
when $J_*\in\cal{J}$.
\end{theo}

\begin{theo}
\n\n
Assume that
 \begin{description} 
\item[$(i)$] $J_* \in \cal{J},$
\item[$(ii)$] $p\rightarrow\infty,~n\rightarrow\infty,~\frac{p}{n}\rightarrow c\in(0,1),$
\item[$(iii)$] If ${}^{\forall}J\in\cal{J}_-$ then $\Omega_J=np\Xi_J=O_p(np)$,
$\lim_{\frac{p}{n}\rightarrow c}\Xi_J=\Xi_{J}^*$,  and $\Xi_{J}^*$ 
is positive definite, 
\item[$(iv)$] ${}^{\forall}J\in\cal{J}$C${}^{\exists}i,j\in J_*,~~
\lim_{\frac{p}{n}\rightarrow c}\ten{\tjic}\inv{\Sigma}\tjic
=\infty,$
 \end{description} 
then
\begin{eqnarray*}
\lim_{\frac{p}{n}\rightarrow c}
\mathrm{P}\left(\argmin_{J\in{\cal J}}\zmcp(J)=J_*\right)=1.
\end{eqnarray*}
\end{theo}

The fact of Theorem 4 is obvious by Section 4 of \cite{davies2006estimation} 
because $\zmcp$ is described by complete sufficient statistics. 
The theorem means that $\zmcp$ is the best unbiased estimator 
of the squared risk of the generalized ridge estimator. 

Theorem 5 is seen as an extension of \cite{fujikoshi2014consistency}. He showed 
that $\mcp$ has consistency under similar conditions to Theorem 5.2 of his paper. 
The difference between the conditions of our result and his is fourth condition.
The assumption $(iv)$ is to ensure that the regression coefficients do not 
have strong correlation, and hence, we can distinguish the candidate models. 
His conditions do not contain the assumption $(iv)$. 
The assumption $(iv)$ is necessary for the consistency of $\zmcp$ 
because $\zmcp$ cannot select the true model when the true model has 
strong correlation of regression coefficients.

{\sc Proof.}\m
From \cite{fujikoshi2014consistency}, we can express the differences between 
$\mcp(J)$ and $\mcp(J_*)$ as
\begin{eqnarray*}
&&\mcp(J)-\mcp(J_*)\\
&=&\left(1-\frac{p+1}{n-k}\right)\left(
(n-k)\left(\tr\left(\lj\inv \mj\right) - \tr\left(\ljt\inv\mjt\right)\right)+2p(\kj-\kjt)\right)\\
&&+p(p+1)\left(\frac{2(\kj-\kjt)}{n-k}\right),
\end{eqnarray*}
where  
\begin{eqnarray*}
\lj&\sim&{\cal W}_{\rj}(p,\I{\rj};\Omega_J),\\
\mj&\sim&{\cal W}_{\rj}(n-\kj-p,\I{\rj}),
\end{eqnarray*}
${\cal W}_p(k,\Sigma;\Omega)$ is the noncentral Wishart distribution that 
has degree of freedom $k$, scale matrix $\Sigma$ and 
noncentral matrix $\Omega$.
Moreover, $\lj$ and $\mj$ are independently distributed for a fixed model $J$
(but $\lj$ and $L_{J'}$ (or $\mj$ and $M_{J'}$) for different $J,~J'\in{\cal J}$ could be depended.).
Based on a well-known asymptotic 
method on Wishart distributions, we can see that under assumptions $(ii)$ and $(iii)$
\begin{eqnarray}
\lim_{\frac{p}{n}\rightarrow c}\frac{1}{np}\lj=\Xi_{J}^*, ~~\lim_{\frac{p}{n}\rightarrow c}
\frac{1}{n}\mj = (1-c)\I{\rj}.
\end{eqnarray}

In the case of $J_* \subset J$, from $(3)$, 
\begin{eqnarray*}
\lim_{\frac{p}{n}\rightarrow c}\frac{1}{n}\lj=c\I{\rj}.
\end{eqnarray*}
Therefore, from $(4)$,
\begin{eqnarray*}
&&\lim_{\frac{p}{n}\rightarrow c}\frac{1}{n}\left\{
\mcp(J)-\mcp(J_*)\right\}\\
&=&(1-c)\left(\frac{c}{1-c}+2c\right)(\kj-\kjt)+2c^2(\kj-\kjt)\\
&=&c(\kj-\kjt)>0
\end{eqnarray*}
in probability.

Similarly, in the case of $J_* \not\subset J$, from $(3)$,
\begin{eqnarray*}
&&\lim_{\frac{p}{n}\rightarrow c}\frac{1}{np}\left\{
\mcp(J)-\mcp(J_*)\right\}\\
&=&(1-c)\frac{1}{1-c}\tr\Xi^*_J=\tr \Xi^*_J >0
\end{eqnarray*}
in probability.

Therefore, it is sufficient to show that 
\begin{eqnarray*}
\lim_{\frac{p}{n}\rightarrow c}\frac{p-2}{n}\tr(\Lj\inv{\Cj})=0,~{}^\forall J\in {\cal J}
\end{eqnarray*}
in probability because
\begin{eqnarray*}
&&\zmcp(J)-\zmcp(J_*)\\
&=&\mcp(J)-\mcp(J_*)+(p-2)\left\{\tr(\Lambda_{J_*}\inv{C_{J_*}})
-\tr(\Lj\inv{\Cj})\right\}.
\end{eqnarray*}
From the definitions of $\xjic$ and $\zf$, 
letting $\eji=(\Dj\ten{\Pj}\ten{\aj}\ajt\Pjt\tjt)_{i,:}$, 
where $A_{i,:}$ is the $i$-th row of $A$, 
we have  
\begin{eqnarray*}
&\dji^{-\frac{1}{2}}\Sigma_{}^{-\frac{1}{2}}\xjic\sim N_T(\dji^{-\frac{1}{2}}\Sigma^{-\frac{1}{2}}\eji,\mathrm{I}),~~~~~
\zf\Sigma_{}^{-\frac{1}{2}}\sim N_{(n-\kf)\times T}(0,\mathrm{I}\otimes\mathrm{I}),\\
&\Sigma_{}^{-\frac{1}{2}}S_F\Sigma_{}^{-\frac{1}{2}}\sim W_T(n-\kf,\mathrm{I}).
\end{eqnarray*}
Therefore we can bound the magnitude of $\lji\inv{\cji}$ as 
\begin{eqnarray*}
&&\lji\cji^{-1}\\
&=&\frac{\dji(p-2)}{n-\kj-p+3}(\ten{\xjic}\isf\xjic)^{-1}\\
&=&\frac{p-2}{n-\kj-p+3}\frac{n-\kf-p-1}{p}
\inv{\left(\frac{\dji^{-1}(n-\kf-p-1)}{p}\ten{\xjic}\isf\xjic\right)}\\
&=&\frac{p-2}{n-\kj-p+3}\frac{n-\kf-p-1}{p}
\frac{\chi^2_{(n-\kf-p-1)}}{n-\kf-p-1}
\frac{p}{\chi^2_p(\dji^{-1}\ten{\eji}\inv{\Sigma}\eji)}\\
&=&\frac{p-2}{n-\kj-p+3}\frac{n-\kf-p-1}{\dji^{-1}\ten{\eji}\inv{\Sigma}\eji+p}
\frac{\chi^2_{(n-\kf-p-1)}}{n-\kf-p-1}
\frac{p+\dji^{-1}\ten{\eji}\inv{\Sigma}\eji}{\chi^2_p(\dji^{-1}\ten{\eji}\inv{\Sigma}\eji)}\\
&=&O_p\left(\frac{p}{\dji^{-1}\ten{\eji}\inv{\Sigma}\eji+p}\right)\\
&=&O_p\left(\frac{p}{n\ten{\eji}\inv{\Sigma}\eji+p}\right)
\end{eqnarray*}
because
\begin{eqnarray*}
&&\inv{\left(\frac{\dji^{-1}(n-\kf-p-1)}{T}\ten{\xjic}\isf\xjic\right)}\\
&\sim& F''(n-\kf-T-1,T,0,\dji^{-1}\ten{\eji}\inv{\Sigma}\eji),
\end{eqnarray*}
where $\chi^2_k(\delta)$ is the noncentral chi-squared distribution that has 
degree of freedom $k$ and non-central parameter $\delta$, in particular, 
we simply write $\chi^2_k$ for $\chi^2_k(0)$ and call $\chi^2_k$ the chi-squared distribution, 
and $F''(n_1,n_2,\gamma_1,\gamma_2)$ is the doubly noncentral $F$ distribution that 
has degree of freedom $(n_1,n_2)$ and non-central parameters $(\gamma_1,\gamma_2)$.
From the assumption $(iv)$ and 
$\ten{\eji}\inv{\Sigma}\eji=
O(1)\times\left(
\sum_{i=1}^{\kj}\sum_{j=1}^{\kj}\ten{\tjic}\inv{\Sigma}\tjic\right)$,
\begin{eqnarray*}
\lim_{\frac{p}{n}\rightarrow c}\frac{p-2}{n}\tr(\Lj\inv{\Cj})=
\lim_{\frac{p}{n}\rightarrow c}\frac{p-2}{n}\sum_{i=1}^{\kj}\lji\cji^{-1}=0
\end{eqnarray*}
in probability.\qed

\subsection{Modified AICc}
\m\n In this subsection, we propose a model selection criterion that is a 
modification of $\aicc$ under the Kullback-Leibler risk. $\aicc$ is the estimator 
of the Kullback-Leibler risk with the maximum likelihood estimator and given by 
\begin{eqnarray*}
\aicc(J)=n\log\left|\frac{1}{n}S_J\right|+np\log2\pi+\frac{np(n+\kj)}{n-\kj-p-1},
\end{eqnarray*}
where $\aicc(J)$ is $\aicc$ under a model $J$.
As in $\mcp$, 
we consider the generalized ridge estimator and 
construct an unbiased estimator of the Kullback-Leibler risk corresponding to that.
$\klrisk{\bjt}{\Sigma}{\hat{f}(\cdot|\bjbay,\Sigma_J)}$ is given by
\begin{eqnarray*}
&&2\klrisk{\bjt}{\Sigma}{\hat{f}(\cdot|\bjbay,\Sigma_J)}\\
&=&2\klrisk{\bjt}{\Sigma}{\hat{f}(\cdot|\bjmle,\Sigma_J)}
-n\frac{p-2}{n-\kj-p+1}\E{\tr(\Lj\inv{\Cj})}
\end{eqnarray*}
by the proof of Theorem 3. 
In contrast to $\aicc$ which is the 
unbiased estimator of 
$\klrisk{\bjt}{\Sigma}{\hat{f}(\cdot|\bjmle,\Sigma_J)}$,
we denote by $\zklic$ an unbiased estimator of 
$\klrisk{\bjt}{\Sigma}{\hat{f}(\cdot|\bjbay,\Sigma_J)}$ 
under a model $J$ which is given by
\begin{eqnarray*}
&&\zklic(J)\\
&=&\aicc(J)-n\frac{p-2}{n-\kj-p+1}\tr(\Lj\inv{\Cj})\\
&=&n\log\left|\frac{1}{n}S_J\right|+np\log2\pi+\frac{np(n+\kj)}{n-\kj-p-1}
-n\frac{p-2}{n-\kj-p+1}\tr(\Lj\inv{\Cj}).
\end{eqnarray*}

The model selection criterion has the following properties.
\begin{theo}
\n\n
$\zklic$ is a uniformly minimum variance unbiased estimator 
of the Kullback-Leibler risk when 
$J_*\subset J$.
\end{theo}

\begin{theo}
\n\n
Assume that 
 \begin{description} 
\item[$(i)$] $J_* \in \cal{J},$
\item[$(ii)$] $p\rightarrow\infty,~n\rightarrow\infty,~\frac{p}{n}\rightarrow c\in(0,1),$
\item[$(iii)$] If $J_*\not\subset J$ then $\Omega_J=np\Xi_J=O_p(np)$,
$\lim_{\frac{p}{n}\rightarrow c}\Xi_J=\Xi_{J}^*$ and $\Xi_{J}^*$
is positive definite,
\item[$(iv)$] ${}^{\forall}J\in{\cal J}$C${}^{\exists}i,j\in J_*,~~
\lim_{\frac{p}{n}\rightarrow c}\ten{\tjic}\inv{\Sigma}\tjic
=\infty,$
 \end{description} 
then
\begin{eqnarray*}
\lim_{\frac{p}{n}\rightarrow c}
\mathrm{P}\left(\argmin_{J\in{\cal J}}\zklic(J)=J_*\right)=1.
\end{eqnarray*}
\end{theo}

The assumption $(iv)$ is made for the same reason as 
discussed in Theorem 5. We again observe that $\zklic$ has 
consistency like $\zmcp$.

{\sc Proof.}\m
From the proof of Theorem 2 of \cite{fujikoshi2012high}, 
\begin{eqnarray*}
\lim_{\frac{p}{n}\rightarrow c}
\frac{1}{n\log{p}}\left\{\aicc(J)-\aicc(J_*)\right\}= \gamma_J>0
,~ J_* \not\subset J,\\
\lim_{\frac{p}{n}\rightarrow c}
\frac{1}{p}\left\{\aicc(J)-\aicc(J_*)\right\}=\rj\left\{\frac{1}{c}\log(1-c)+2\right\}>0
,~ J_* \subset J,~J\not=J_*.
\end{eqnarray*}
Therefore, it is sufficient to show that
\begin{eqnarray}
\lim_{\frac{p}{n}\rightarrow c}\frac{1}{n\log p}
\frac{n(p-2)}{n-\kj-T+1}\tr(\Lj\inv{\Cj})=0,
J_* \not\subset J,\nonumber\\
\lim_{\frac{p}{n}\rightarrow c}\frac{1}{p}\frac{n(p-2)}{n-\kj-T+1}\tr(\Lj\inv{\Cj})=0,
~J_* \subset J,~J\not=J_*
\end{eqnarray}
in probability because 
\begin{eqnarray*}
\zklic=\aicc-\frac{n(p-2)}{n-\kj-p+1}\tr(\Lj\inv{\Cj}).
\end{eqnarray*}
In the case of $J_*\subset J$, $J\not=J_*$, 
$(5)$ can be easily shown by the assumption $(iv)$ and 
the proof of Theorem 5.

Let $J_* \not\subset J$ then
\begin{eqnarray*}
\inv{\left(\frac{\dji^{-1}(n-\kj-p-1)}{p}\ten{\xjic}\isj\xjic\right)}
\sim F''(n-\kj-p-1,p,0,\dji^{-1}\ten{\eji}\inv{\Sigma}\eji)
\end{eqnarray*}
because
\begin{eqnarray*}
&\dji^{-\frac{1}{2}}\Sigma_{}^{-\frac{1}{2}}\xjic\sim N_p(\dji^{-\frac{1}{2}}\Sigma^{-\frac{1}{2}}\eji,\mathrm{I}),~~~~~
\Sigma_{}^{-\frac{1}{2}}S_J\Sigma_{}^{-\frac{1}{2}}\sim W_p(n-\kj,\mathrm{I},\Omega_J).
\end{eqnarray*}
Therefore, from the assumption $(iv)$ 
\begin{eqnarray*}
\lim_{\frac{p}{n}\rightarrow c} \frac{n}{n\log p}
\frac{p-2}{n-\kj-p+1}\tr(\Lj\inv{\Cj})=0
\end{eqnarray*}
in probability. 
\qed\\

\section{Numerical study}
\m\n In this section, we numerically examine the validity of our propositions. 
The risk of a selected model and the probability of selecting the true model
by $\mcp$, $\zmcp$, $\aic$, $\aicc$, and $\zklic$ were 
evaluated by Monte Carlo simulations with 1,000 iterations. 
The ten candidate models $J_\alpha = \{1,...,\alpha\}~(\alpha=1,...,10)$ were
evaluated. 
In the experiment to evaluate the risks of the selected models, we employed 
$n=100,200$ and $p/n=0.04,0.06,...,0.8$.
In the experiment to evaluate the probability of selecting the true 
model, we employed $n=100,200,400,600$ and $p/n=0.04,0.06,...,0.8$.
The true model was determined by $\bjt=\ten{(1,-2,3,-4,5)}\ten{1_p}$, 
$J_*=\{1,2,3,4,5\}$, and 
the $(a,b)$-th element of $\Sigma$ was defined by 
$(0.8)^{|a-b|}~(a=1,...,p;b=1,...,p)$.
Here, $1_p$ is the $p$-dimensional vector of ones. 
Thus, $J_1,J_2,J_3$, and $J_4$ are underspecified models, and 
$J_6,J_7,J_8,J_9$, and $J_{10}$ are overspecified models. 
Explanatory variables $A$ is generated in two different ways: 
in Case 1, $A_{a,b} = u_{a}^{(b-1)}~(a=1,...,n;b=1,...,10)$, where 
$u_1,\dots,u_n \sim U(-1,1)$ i.i.d, and in Case 2, $A_{a,b}
\sim N(0,1)$ i.i.d..

\renewcommand{\figurename}{Figure\n}

\begin{figure}[H]
  \begin{center}
   \includegraphics[width=80mm]{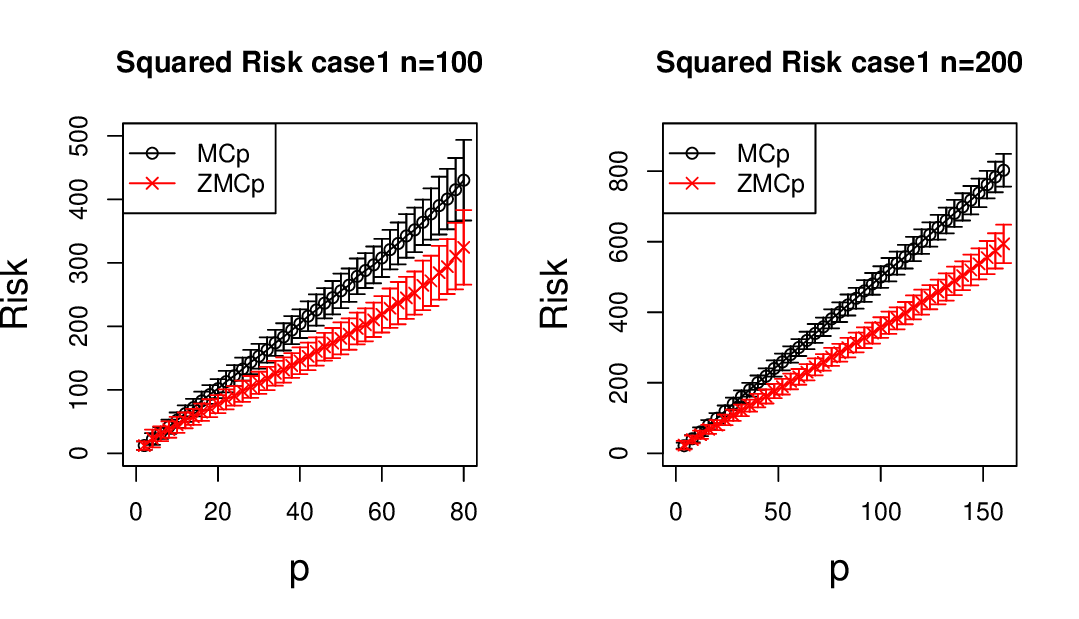}
  \caption{Comparison between $\mcp$ and $\zmcp$ 
  under the squared risk in 
 \newline \n\n\n\n\n\n\n\n\n\n case 1.}
  \label{fig:one}
  \end{center}
\end{figure}

\begin{figure}[H]
  \begin{center}
   \includegraphics[width=80mm]{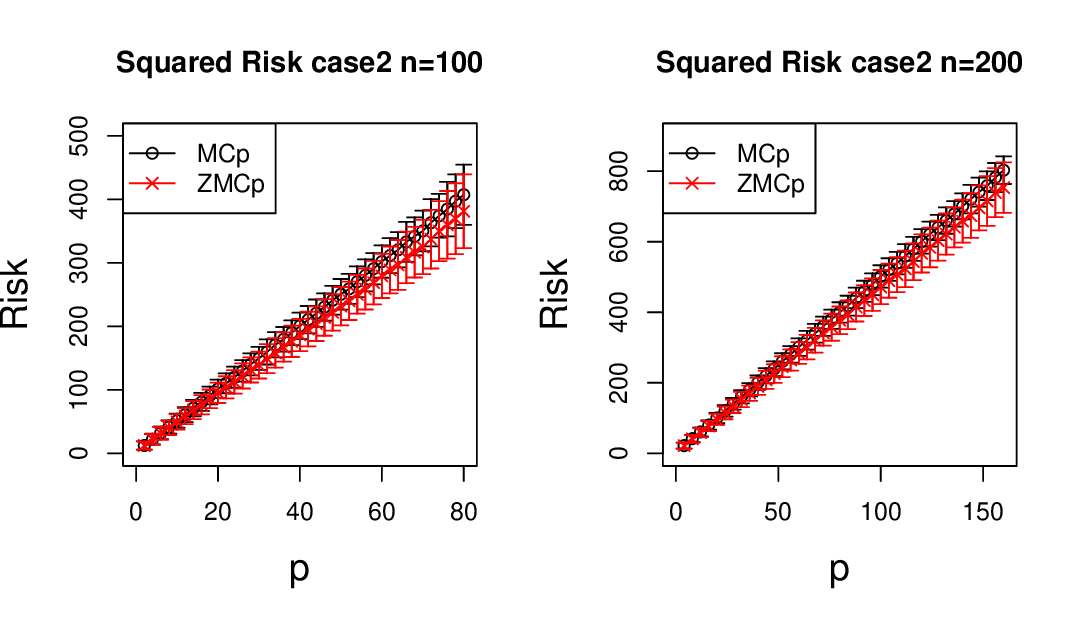}
  \caption{Comparison between $\mcp$ and $\zmcp$ 
   under the squared risk in 
\newline \n\n\n\n\n\n\n\n\n\n case 2.}
  \label{fig:two}
  \end{center}
\end{figure}

Figure 1 and Figure 2 show the squared risks of the selected models by 
$\mcp$ and $\zmcp$. 
While Figure 3 and Figure 4 show the Kullback-Leibler risks 
of the selected models by $\aic$, $\aicc$, and $\zklic$,
the constant part depending on the true distribution was subtracted.
In Case 1, it is seen that $\zmcp$ largely improves $\mcp$. 
On the other hand, in Case 2, the difference between the squared risks of 
the selected models is not as much as that in Case 1. 
The reason is that 
the explanatory variables have larger correlation in Case 1 than in Case 2, 
and the generalized ridge estimator is more robust against the correlation.  

\begin{figure}[H]
  \begin{center}
   \includegraphics[width=80mm]{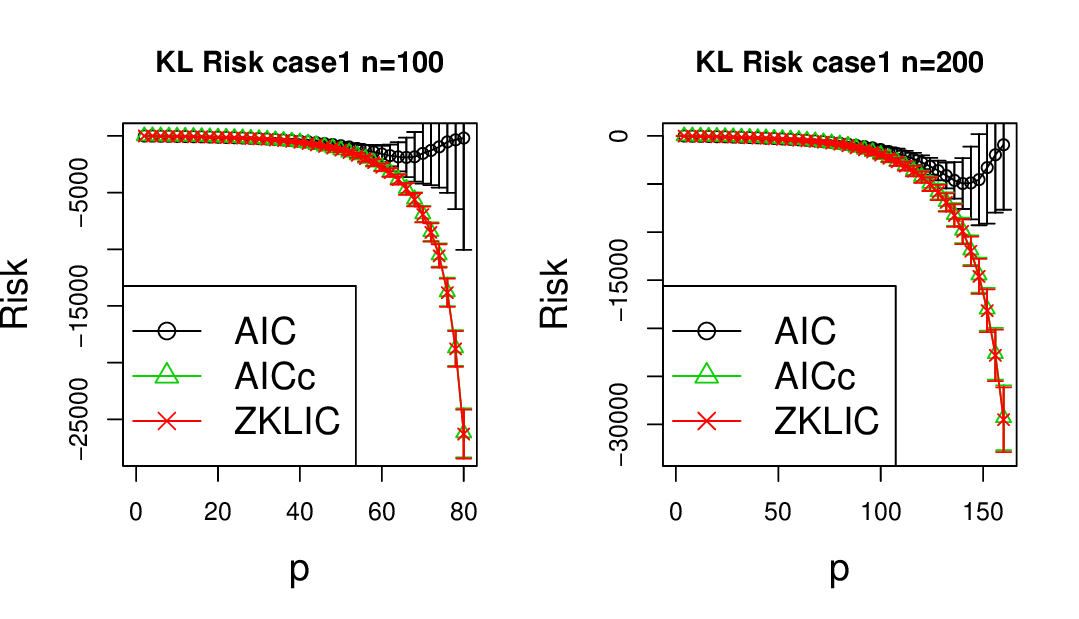}
  \caption{Comparison between $\aic$, $\aicc$ and   
   $\zklic$ under the   
  \newline \n\n\n\n\n\n\n\n\n\n Kullback-Leibler risk in case 1.}
  \label{fig:one}
  \end{center} 
\end{figure}

\begin{figure}[H]
  \begin{center}
   \includegraphics[width=80mm]{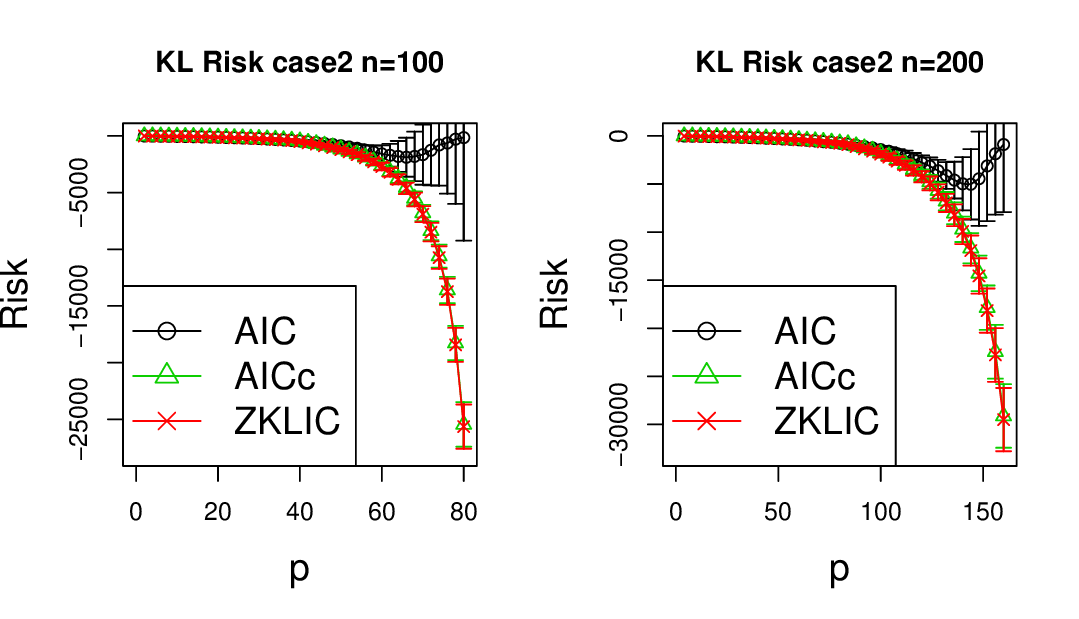}
  \caption{Comparison between $\aic$, $\aicc$ and   
$\zklic$ under the   
  \newline \n\n\n\n\n\n\n\n\n\n Kullback-Leibler risk in case 2.}
  \label{fig:two}
  \end{center}
\end{figure}

Figure 5 and Figure 6 show the probability of selecting the true model 
by each model selection criterion. 
In Figure 6, the probability of selecting the true model by our model selection criteria 
is large when the sample size is large. 
However, in Figure 5, this probability is small; the matrix of regression coefficient 
has large correlation. 
Furthermore, the probabilities of selecting the true model by 
our model selection criteria are smaller than those of existing ones in 
each case. The reason for this is that the variance of our model selection 
criteria is bigger than those of the existing ones.

\begin{figure}[H]
  \begin{center}
   \includegraphics[width=100mm]{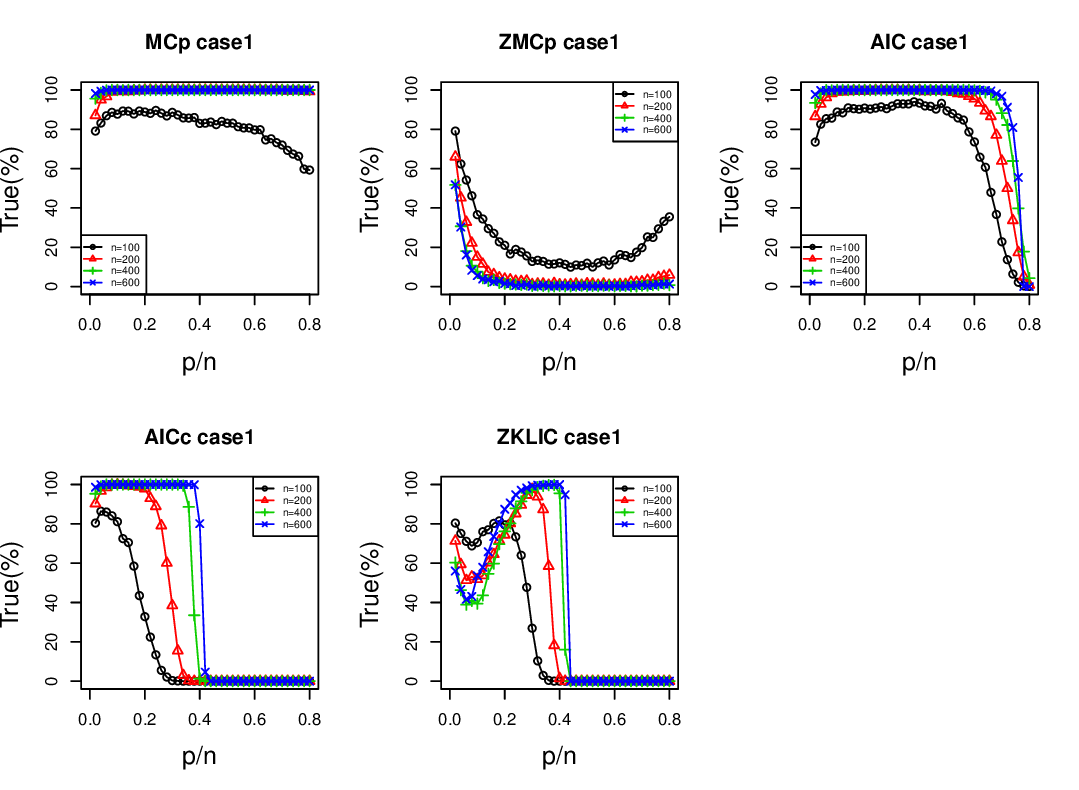}
\caption{Comparison between $\mcp,~\zmcp,~\aic,~\aicc,$ and $\zklic$ of 
\newline \n\n\n\n\n\n\n\n\n\n the probability of selecting the true model in case 1}
  \end{center}
\end{figure}

\begin{figure}[H]
  \begin{center}
   \includegraphics[width=100mm]{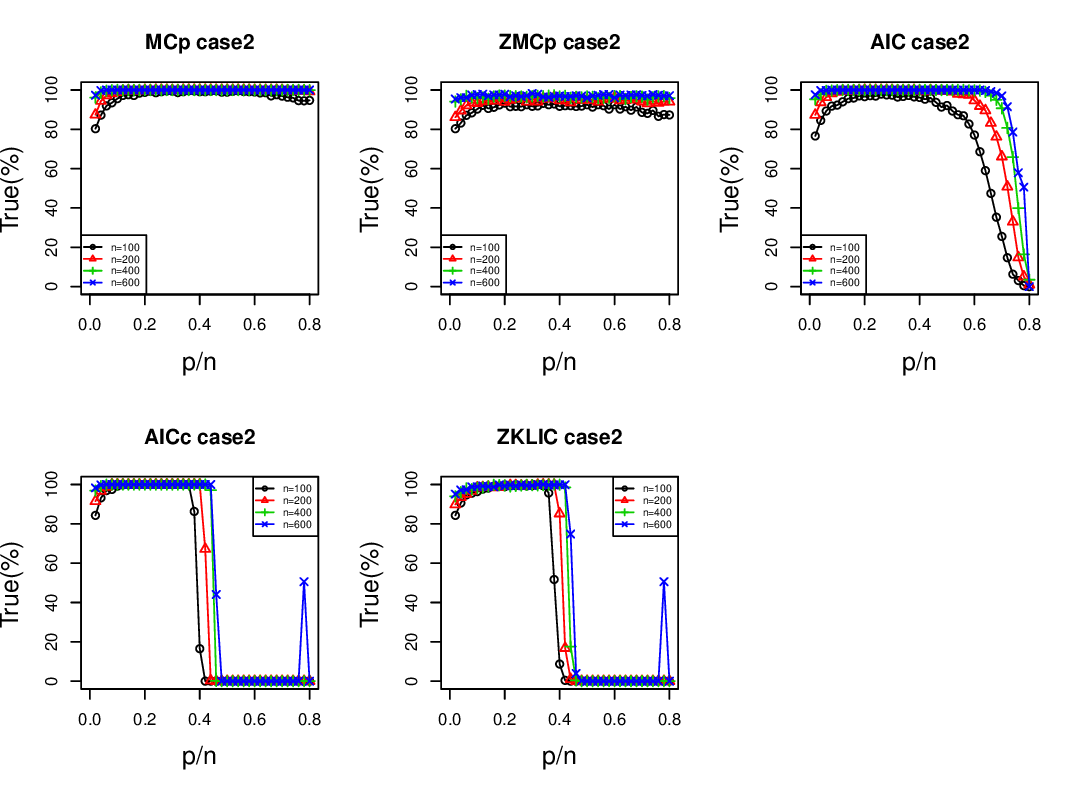}
\caption{Comparison between $\mcp,~\zmcp,~\aic,~\aicc,$ and $\zklic$ of 
\newline \n\n\n\n\n\n\n\n\n\n the probability of selecting the true model in case 2}
  \end{center}
\end{figure}

Although the risks of our model selection criteria are smaller than the ones based on the 
maximum likelihood estimator, the probability of selecting the true model 
by our criteria is worse than their probabilities. This is because 
predictive efficiency and consistency are not compatible (\cite{yang2005can}) 
and our criteria are specialized in making the 
risk smaller. 

\section{Conclusion}
\m\n In this paper, we proposed model selection criteria based on the generalized 
ridge estimator, which improves the maximum likelihood estimator 
under the squared risk and the Kullback-Leibler risk, 
in multivariate linear regression. Moreover, we showed that 
our model selection criteria have the same properties as 
$\mcp$, $\aic$, and $\aicc$ in a high-dimensional asymptotic framework. 
We demonstrated through the numerical experiments that 
our model selection criteria have better performances in terms of 
the risks than the ones based on the maximum likelihood estimators, 
especially when the matrix of regression coefficients has strong correlation.

\section*{Acknowledgement}
\m\n This work was partially supported by MEXT kakenhi (25730013, 25120012, 
26280009, 15H01678 and  15H05707), JST-PRESTO and JST-CREST.

\end{document}